\documentclass{article}
\usepackage{latexsym}
\usepackage{amssymb}
\usepackage{amsfonts}                     
\usepackage{amsmath}
\usepackage[all]{xy}
\usepackage{color}
\usepackage{mathrsfs}
\usepackage[utf8]{inputenc}
\usepackage[T1]{fontenc}
\usepackage{amscd}
\usepackage{amsxtra}
\usepackage{tikz-cd}
\newtheorem{thm}{\bf Theorem}[section]
\newtheorem{cor}[thm]{\bf Corollary}
\newtheorem{lem}[thm]{\bf Lemma}
\newtheorem{prop}[thm]{\bf Proposition}
\newtheorem{defn}[thm]{\bf Definition}
\newtheorem{rem}[thm]{\bf Remark}

\newtheorem{exmp}[thm]{\bf Example}

\newcommand{\field}[1]{\mathbb{#1}}

\newcommand{\N }{\field{N}}
 
\usepackage{tikz}
\usetikzlibrary{positioning}

\def\E{{\cal E}}
\def\F{{\cal F}}

\def\Ext{{\rm Ext}}
\def\Hom{{\rm Hom}}

\def\proof{{\parindent0pt {\bf Proof.\ }}}

\def\pd{{\rm pd}}
\def\fd{{\rm fd}}

\def\wdim{{\rm wdim}}
\def\gldim{{\rm gldim}}

\newcommand{\cqfd}
{\hspace{1cm}
\rule{2mm}{2mm}%
\medbreak%
\par%
}

\begin{document}

\title{When every $S$-flat module is (flat) projective}
\author{Driss Bennis and Ayoub Bouziri}

\date{}
\maketitle
\maketitle
\begin{abstract} Let $R$ be a commutative ring with identity and $S$ a multiplicative subset of $R$. The aim of this paper is to study the class of commutative rings in which every $S$-flat module is flat (resp.,  projective). An $R$-module $M$ is said to be $S$-flat if the localization of $M$ at $S$, $M_S$, is a flat $R_S$-module. Commutative rings $R$ for which all $S$-flat $R$-modules are flat are characterized by the fact that $R/Rs$ is a von Neumann regular ring for every $s\in S$. While, commutative rings $R$ for which all $S$-flat $R$-modules are projective are characterized by the following two conditions: $R$ is perfect and the Jacobson radical $J(R)$ of $R$ is $S$-divisible. Rings satisfying these conditions are called $S$-perfect. Moreover, we give some examples to distinguish perfect rings, S-perfect rings, and semisimple rings. We also investigate the
transfer results of the "S-perfectness" for various ring constructions, which allows the construction of more interesting examples.
\end{abstract}
\medskip



\section{Introduction}
 \hskip .5cm  Throughout this paper, $R$ is a commutative ring with identity, all modules are unitary and $S$ is a multiplicative subset of $R$; that is,  $1 \in S$ and $s_1s_2 \in S$ for any $s_1,s_2 \in S$. Unless explicitly stated otherwise, when we consider a multiplicative subset $S$ of $R$, then we implicitly suppose that $0\notin S$. This will be used in the sequel without explicit mention. We denote the set of all regular elements, von Neumann regular elements and units in $R$ by $reg(R)$, $vnr(R)$ and $U(R)$, respectively. An element $a \in R$ is said to be von Neumann regular if $a^2x = a$ for some $x \in R$ \cite{Abu1, And1}. Recall that $reg(R) = \{a \in R : ann_R(a) = 0\}$, where $ann_R(a)$ is the annihilator of an element $a \in R$, $vnr^*(R)=vnr(R)\setminus \{0\}$ and $U(R)$ are examples of multiplicative subsets. Let $M$ be an $R$-module. As usual, we use $\pd_R(M)$ and $\fd_R(M)$ to denote, respectively, the classical projective dimension and flat dimension of $M$, and, $\wdim(R)$ and $\gldim(R)$ to denote, respectively, the weak and global homological dimensions of $R$. Also, we denote by $M_S$ the localization of $M$ at $S$. Recall that $M_S \cong M \otimes_R R_S$.

\medskip 
In the last years, the notion of $S$-property draw attention of several authors. This  notion  was introduced in 2002 by D. D. Anderson and Dumitrescu where they defined the notions of $S$-finite modules and $S$-Noetherian rings. Namely, an $R$-module $M$ is said to be $S$-finite module if there exist a finitely generated submodule $N$ of $M$ and $s \in S$ such that $sM \subseteq N$. A commutative ring $R$ is said to be $S$-Noetherian if every ideal of $R$ is $S$-finite. In \cite{Benn1}, Bennis and El Hajoui investigated an $S$-versions of finitely presented modules and coherent rings which are called, respectively, $S$-finitely presented modules and $S$-coherent rings. An $R$-module $M$ is said to be $S$-finitely presented if there exists an exact sequence of $R$-modules $0 \to K\to L\to  M \to 0$, where $L$ is a  finitely generated free $R$-module and $K$ is an $S$-finite R-module. Moreover, a
 commutative ring $R$ is called $S$-coherent, if every finitely generated ideal of $R$ is $S$-finitely presented. They showed that the $S$-coherent rings have a similar characterization to the classical one given by Chase for coherent rings  \cite[Theorem 3.8]{Cha1}. Subsequently, they asked whether there exists an $S$-version of Chase's theorem \cite[Theorem 2.1]{Cha1}. In other words, how to define an $S$-version of flatness that characterizes $S$-coherent rings similarly to the classical case?  This problem was solved by the  notion    of $S$-flat module in \cite{Qi1}. Recall that an $R$-module $M$ is said to be $S$-flat if for any finitely generated ideal $I$ of $R$, the natural homomorphism $(I \otimes_{R} M)_S \to (R \otimes_{R} M)_S$ is a monomorphism \cite[Definition 2.5.]{Qi1}; equivalently $M_{S}$ is a flat $R_{S}$-module \cite[Proposition 2.6]{Qi1}. Notice that any flat $R$-module is $S$-flat.

\medskip
 In 2019, Bazzoni and Positselski \cite{Baz1} introduced the notion of an $S$-almost perfect ring $R$; that is,  $R_S$ is a perfect ring and $R/sR$ is a perfect ring for every $s \in S$ \cite[Definition 7.6.]{Baz1}. They showed that a ring $R$ is $S$-almost perfect if and only if every flat $R$-module is $S$-strongly flat. An $R$-module $F$ is said to be $S$-strongly flat if $\Ext^1_R(F, C) = 0$ for all $S$-weakly cotorsion $C$ (i.e., $\Ext^1_R(R_S, C) = 0$). A general framework for $S$-strongly flat modules was developed in the paper \cite{Pos1}. Obviously, we have the following implications:

\begin{center}
projective $\Rightarrow$ $S$-strongly flat $\Rightarrow$ flat $\Rightarrow$ $S$-flat
\end{center}

\medskip 
 In this paper, we study the class of rings in which every $S$-flat $R$-module is flat (resp.,  projective). We define and discuss $S$-weakly\footnote{One can see that the fact that every $R$-module is $S$-flat is equivalent to say that $R_S$ is a von Neumann regular ring. So there is no reason to call such a ring $S$-von Neumann regular.} von Neumann regular rings, which form a wider class than that of von Neumann regular rings, as rings for which every $S$-flat $R$-module is flat. We prove that a ring $R$ is $S$-weakly von Neumann regular if and only if $R/Rs$ is von Neumann regular for every $s\in S$ (Theorem \ref{2thm-S-VNR}). A commutative ring $R$ for which all $S$-flat $R$-modules are  projective are characterized by the following two conditions: $R$ is perfect and the Jacobson radical $J(R)$ of $R$ is $S$-divisible (Theorem \ref{3thm-S-perfect-jacobson}). Recall that an $R$-module $D$ is $S$-divisible if $sD =D$ for every $s \in S$. Rings satisfying these conditions are called $S$-perfect. These rings form a wider class than that of semisimple rings and a smaller class than that of perfect rings.
      
\medskip 
The organization of the paper is as follows: In Section 2, some elementary
properties of $S$-flat modules are obtained. We introduce and study an $S$-version of projective modules. We call it an $S$-projective module (Definition \ref{2-def-S-proj}). Then, we study the behaviour of $S$-flatness and $S$-projectivity on direct product of rings (see Corollaries \ref{2cor-S-projandS-flatinproduct} and \ref{2cor-S-projandS-flatinproduct-pratique}). In section 3, we study when every $S$-flat $R$-module is flat. Rings satisfying this condition are called $S$-weakly von Neumann regular (Definition \ref{2def-S-VNR}). The main result of this section is Theorem \ref{2thm-S-VNR}. In Section 4, we deal with  $S$-perfect rings. We shows that a ring $R$ is $S$-perfect if and only if $R$ is $S$-almost perfect and $S$-von Neumann regular, and $R_S$ is a projective $R$-module (see Theorem \ref{3thm-S-perf-almostperf-vnr}). In Theorem \ref{3thm-S-perfect-jacobson}, we prove that a ring $R$ is $S$-perfect if and only if, $R$ is perfect and the Jacobson radical $J(R)$ of $R$ is $S$-divisible. Finally, we investigate the transfer of the "$S$-perfect" property along a ring homomorphism (Proposition \ref{3prop-S-pefect-in-homomorphism}), and we examine conditions under which $R\ltimes M$  is an $S$-perfect ring, where $S$ is a multiplicative subset of $R\ltimes M$.  We found that, unlike the "perfect" property \cite[Proposition 1.15., page 23]{Fos1}, the "$S$-perfect" property depends on the structure of the $R$-module $M$ considered in the trivial extension $R\ltimes M$ (see Theorem \ref{3thm-s-perfect-in-trival-extention} and Corollary \ref{3cor-s-perfect-in-trival-extention}).   

\section{S-flat and S-projective modules}
 Recall that an $R$-module $M$ is said to be $S$-flat if for any finitely generated ideal $I$ of $R$, the natural homomorphism $(I \otimes_{R} M)_S \to (R \otimes_{R} M)_S $ is a monomorphism; equivalently, $M_{S}$ is flat $R_{S}$-module \cite[Proposition 2.6]{Qi1}. We denoted by $S\F$ the class of all $S$-flat modules. It is well-know that the class $\F$ of flat modules is closed under pure submodules, pure quotient modules, extensions and directed limits. Here we have the corresponding result for $S$-flat modules.

\begin{lem}\label{2-recall-vnr}
$S\F$ is closed under extensions, direct sums, direct summands, direct limits pure submodules and pure quotients.
\end{lem}

\proof Recall that an $R$-module $M$ is $S$-flat if and only if $M_S\cong R_S\otimes_R M$ is a flat $R_S$-module \cite[Proposition 2.6]{Qi1}. Thus, all properties follow from their validity for the class of flat $R_S$-modules and the fact that they are preserved by the exact functor $R_S\otimes_R(-)$.  \cqfd
 The canonical ring homomorphism $\theta: R \to R_{S}$ makes every $R_{S}$-module an $R$-module via the module action $r.m = \frac{r}{1}.m$, where $r\in R $ and $m\in M$. Recall from \cite{Rot1} the following lemma that we frequently use in this paper.
     
\begin{lem}[\cite{Rot1}, Corollary 4.79]\label{2-lem-localiz-isomorphism}
Every $R_S$-module $M$ is naturally isomorphic to its localization $M_S$ as $R_S$-modules \end{lem}


In order to facilitate the use of $S$-flat modules in next sections, we propose the following characterisation:
\begin{prop}\label{2-prop-S-flat-charac} An $R$-module $M$ is $S$-flat if and only if  $M_S$ is a flat $R$-module.
\end{prop} 

 \proof Assume that $M$ is an $S$-flat $R$-module, then $M_S$ is a flat $R_S$-module. Since $R_S$ is a flat $R$-module, $M_S$ is a flat $R$-module \cite[Corollary 10.72]{Rot1}.  Conversely, if $M_S$ is a flat $R$-module, then $(M_S)_S$ is a flat $R_S$-module. By Lemma \ref{2-lem-localiz-isomorphism} $(M_S)_S\cong M_S$. Therefore, $M_S$ is a flat $R_S$-module. Hence, $M$ is an $S$-flat $R$-module. \cqfd  

 In this paper, we are also interested in the $R$-projectivity of the localization of an $R$-module $M$ at $S$. This is why we introduce an $S$-version of projectivity.

\begin{defn}\label{2-def-S-proj} An $R$-module $M$ is said to be $S$-projective if $M_S$ is a projective $R$-module.
\end{defn}

If $S$ is composed of units, then $S$-projective modules are exactly projective modules. It follows easily from Proposition \ref{2-prop-S-flat-charac} that every $S$-projective module is $S$-flat.
Combining Propositions \ref{2-prop-S-flat-charac} and Lemma \ref{2-lem-localiz-isomorphism}, one can see easily that an $R_S$-module is $S$-flat as $R$-module if and only if it is a flat $R$-module. Using similar arguments, the corresponding result is also true if we consider $S$-projective modules.
\begin{prop}\label{2prop-S-proj-R_S}  An $R_S$-module $M$ is $S$-projective if and only if it is projective as $R$-module.
\end{prop}

 Recall that an $R$-module $F$ is said to be $S$-strongly flat if $\Ext^1_R(F, C) = 0$ for all $S$-weakly cotorsion $C$ (i.e., $\Ext^1_R(R_S, C) = 0$) \cite{Pos1}. It is clear that any projective $R$-module is $S$-strongly flat. The following result characterizes when the converse holds in term of $S$-projectivity.
\begin{prop}\label{2prop-proj-are-S-proj} The following assertions are equivalent:
\begin{enumerate}
\item $R_S$ is a projective $R$-module. 
\item Every projective $R$-module is $S$-projective.
\item Every $S$-strongly flat is projective.
\end{enumerate}
\end{prop}

 \proof $1.\Rightarrow 2.$ Let $M$ be a projective $R$-module and $\E$ an exact sequence of $R$-module. The exactness of the sequence $ \Hom_R(M_S ,\E)\cong \Hom_R(M\otimes_R R_S ,\E)$ follows form the natural isomorphism $$\Hom_R(M, \Hom_R(R_S,\E)) \cong  \Hom_R(M\otimes_R R_S ,\E)$$ \cite[Theorem 1.1.8]{Gla1} and the assumption on $R_S$.  
Implications $2. \Rightarrow 1.$ and $3.\Rightarrow 1. $ are obvious. And the implication $1. \Rightarrow 3.$ follows by \cite[Lamma 3.1]{Baz1}.\cqfd


Recall that a ring $R$ is $S$-almost perfect if and only if every flat $R$-module is $S$-strongly flat \cite[Theorem 7.9]{Baz1}.
\begin{cor}\label{2cor-perfect-S-almost-perfect}
  The following assertions are equivalent:
\begin{enumerate}
\item $R$ is $S$-almost perfect and $R_S$ is projective.  
\item $R$ is a perfect ring.
\end{enumerate}
\end{cor}
 
 \proof $1.\Rightarrow 2.$ This follows by Proposition \ref{2prop-proj-are-S-proj} and the fact that a ring $R$ is perfect if and only if every flat $R$-module is projective \cite[Theorem 3.10.22]{Wan1}. 
 
 $.2\Rightarrow 1.$  Recall that $R_S$ is a flat $R$-module \cite[Theorem 4.80]{Rot1}. Since $R$ is perfect, $R_S$ is projective by \cite[Theorem 3.10.22]{Wan1}. \cqfd

We end this section by characterizing the $S$-flatness and $S$-projectivity in the direct product of rings. Let $R = R_1 \times\cdots\times R_n$ be the direct product of commutative rings $R_1,...,R_n$. Let $e_i = (0, ..., 1, ..., 0)$; that is,  the $i$-th component of $e_i$ is $1$ and the others are $0$. Let $M$ be an $R$-module. Then $M = M_1\times\cdots \times M_n$, where $M_i=e_iM$ is an $R_i$-module for each $i=1,...,n$ (see \cite[Example 1.8.21]{Wan1}).  It is well-known that $M$ is a projective $R$-module (resp., a flat) if and only if $M_i$ is a projective (resp., a flat) $R_i$-module for each $i=1,...,n$ (see \cite[Exercise 1.12., page 64]{Wan1}). The following result will be needed in the next sections.

\begin{prop}\label{2prop-locali-of-product} Let $R = R_1\times R_2 \times\cdots\times R_n$ be the direct product of commutative rings $R_1,R_2,..., R_n$, and $S$ be a multiplicative set of  $R$. For an $R$-module $M$, there is an isomorphism of $R$-modules: $$M_S \cong e_1M_{e_1S}\times e_2M_{e_2S}\times\cdots\times e_n M_{e_nS},$$ with $e_iM_{e_iS}=0$ if $0\in e_iS$. 

In particular, if $S_i$ is a multiplicative subset of $R_i$ for $i = 1, 2, ..., n$, then  $M_S \cong M_{S_1}\times M_{S_2}\times\cdots\times M_{S_n}$.
\end{prop}

\proof  We may assume that there is $n_0< n$ such that $0\in e_i S \Leftrightarrow i< n_0 $. Let $M=e_1M\times \cdots \times e_n M$ be an $R$-module.  Define 
$$\phi : M_S \to e_{n_0} M_{e_{n_0} S}\times \cdots\times e_{n}M_{e_{n}S}$$ by 
 $\phi(\frac{(m_1,... ,m_{n})}{(s_1,..., s_{n})}) = (\frac{m_{n_0}}{s_{n_0}},..., \frac{m_{n}}{s_{n}})$.  It is routine to check that this is an isomorphism. The last statement follows by taking $S=S_1\times \cdots\times S_n$ and $n_0=1$. \cqfd 
The following result show that $S$-flat modules and $S$-projective modules have the same behaviour as flat modules and projective modules on direct products of rings \cite[Exercise 1.12., page 64]{Wan1}.
\begin{cor}\label{2cor-S-projandS-flatinproduct}
Let $R = R_1 \times R_2\cdots\times R_n$ be the direct product of commutative rings $R_1,...,R_n$, and $S$ be a multiplicative subset of  $R$ and $n_0$ as in the Proof of Proposition \ref{2prop-locali-of-product}. The following hold for an $R$-module $M$:
\begin{enumerate}
\item $M$ is an $S$-flat $R$-module if and only if $e_iM$ is an $e_iS$-flat $R_i$-module for each $i=n_0,..., n.$ 
\item $M$ is an $S$-projective $R$-module if and only if $e_iM$ is an $e_iS$-flat $R_i$-module for each $i=n_0,..., n.$ 
\end{enumerate}
\end{cor}

\proof We only check $1.$,  $2.$ follows similarly. Let $M$ be an $S$-flat $R$-module. Then, by Proposition \ref{2prop-locali-of-product}, $M_S \cong e_{n_0}M_{e_{n_0} S}\times \cdots\times e_{n} M_{e_{n}S}$. Then the result follows by Proposition  \ref{2-prop-S-flat-charac} and \cite[Exercise 1.12., page 64]{Wan1}. \cqfd

We have the following interesting case. 

\begin{cor}\label{2cor-S-projandS-flatinproduct-pratique}
Let $R = R_1 \times R_2\cdots\times R_n$ be the direct product of commutative rings $R_1,...,R_n$, and $S =S_1\times S_2\cdots \times S_n$ be a Cartesian product of multiplicative subsets of $R_1,...,R_2$, respectively. An $R$-module $M =M_1 \times M_2\cdots\times M_n$ is $S$-projective (resp., $S$-flat) if and only if the induced $R_i$-module $M_i$ is $S_i$-projective (resp., $S_i$-flat) for each $i = 1, ...,n$.
\end{cor}

\section{S-weakly von Neumann rings}
In this section we investigate the following notion.
\begin{defn}\label{2def-S-VNR} We say that $R$ is $S$-weakly von Neumann regular if each $S$-flat $R$-module is a flat $R$-module.
\end{defn}

 It is well-know that a ring $R$ is von Neumann regular if and only if every $R$-module is flat \cite[page 212]{Rot1}. Then, every von Neumann regular ring is $S$-weakly von Neumann regular. Moreover, a commutative ring $R$ is von Neumann regular if and only if it is $S$-weakly von Neumann regular and $R_S$ is von Neumann regular. Indeed, this follows from the fact that $R_S$ is von Neumann regular if and only if every $R$-module is $S$-flat (Recall that an $R$-module $M$ is $S$-flat if and only the $R_S$-module $M_S$ is flat). 
 
Before giving the main result of this section which characterizes $S$-weakly von Neumann regular rings, we need to recall some notions. Recall from \cite[Definition 1.6.10]{Wan1} that an $R$-module $M$ is called $S$-torsion if for any $m \in M$, there is an $s\in S$ such that $sm=0$; equivalently $M_S = 0$ (see \cite[Example 1.6.13 ]{Wan1}). Thus, any $S$-torsion module is $S$-flat (see Proposition \ref{2-prop-S-flat-charac}). It is clear that a submodule  of an $S$-torsion module is $S$-torsion.



Recall from \cite[Theorem 1.4.5.]{Gla1} that an element $a\in R$ is von Neumann regular if and only if $R/Ra$ is a flat $ R$-module.
\begin{lem}\label{2-recall-vnr} Let $R$ be a ring and let $a$ be an element of $R$.  If $R/Ra^2$ is a von Neumann regular ring, then $a\in vnr(R)$; equivalently, $R/Ra$ is a flat $ R$-module.
\end{lem}

\proof Since $R/Ra^2$ is a von Neumann regular, there exists $r\in R$ such that $\overline{a}=r\overline{a}^2=0$. So $a\in Ra^2$, as desired. The last assertion follows by \cite[Theorem 1.4.5.]{Gla1}. \cqfd

 The main result of this section is the following:

\begin{thm}\label{2thm-S-VNR}    The following assertions are equivalent: 
\begin{enumerate}
\item $R$ is $S$-weakly von Neumann regular. 
\item Every $S$-torsion module is flat.
\item  $R/Rs$ is a von Neumann regular ring for every $s\in S$. 
\item $S\subseteq vnr^*(R)$ and the kernel of the canonical homomorphism  $\pi_M : M \to M_S$ is flat for any  $R$-module $M$.
\item For every ($S$-flat) $R$-module $M$, the canonical homomorphism $\pi : M \to M_S$ is an epimorphism and its kernel is a flat $R$-module. 
\end{enumerate}
\end{thm}
\proof $1. \Rightarrow 2.$ Let $M$ be an $S$-torsion module. Then $M_S=0$ \cite[Example 1.6.13]{Wan1}, so $M$ is an $S$-flat module. Hence, by 1., $M$ is a flat $R$-module.  

$2. \Rightarrow 3.$ 
Let $M$ be an $R/Rs$-module, then $M$ is an $S$-torsion $R$-module, so $M$ is a flat $R$-module. Then, $M\cong M\otimes_R R/Rs$ is a flat $R/Rs$-module \cite[Theorem 1.2.9]{Gla1}. Hence, $R/Rs$ is a von Neumann regular ring. 

$3. \Rightarrow 2.$ Let $M$ be an $S$-torsion module. By \cite[Proposition 3.48]{Rot1}, it suffices to prove that every finitely generated submodule  $N \subseteq M$ is flat. Since $N$ is $S$-torsion,  $sN=0$ for some $s\in S$ because $N$ is finitely generated. Therefore, $N$ is an $R/Rs$-module. Since $R/Rs$ is von Neumann regular, $N$ is a flat $R/Rs$-module. Since, by Lemma  \ref{2-recall-vnr}, $R/Rs$ is a flat $R$-module, $N$ is a flat $R$-module \cite[Theorem 1.2.8]{Gla1}.


 $2. \Rightarrow 4.$ The first assertion follows by $2. \Rightarrow 3.$ and Lemma \ref{2-recall-vnr}. The second assertion follows from the fact that the kernel of the canonical homomorphism  $\pi_M : M \to M_S$ is $S$-torsion for any  $R$-module $M$.
 
 $4. \Rightarrow 5.$ Let $\frac{m}{s} \in M_S$. Since $s$ is von Neumann regular, there exists $r\in R$ such that $rs^2=s$ 
, so $ms=mrs^2$ and then $\frac{m}{s}=\frac{mr}{1}$. This means that $\pi_M$ is an epimorphism.
 
$5.\Rightarrow 1.$ Let $M$ be an $S$-flat $R$-module.  We have the following exact sequence: 
$$ 0\to K \to M \to M_S\to 0$$

Since $M_S$ is a flat $R$-module by Proposition \ref{2-prop-S-flat-charac} and $K$ is flat by assumption, $M$ is a flat $R$-module. \cqfd
 
  We have the following consequence:  

\begin{cor}\label{2cor-S-VNR}  Suppose $S\subseteq  reg(R)$, then the following assertions are equivalent: 
\begin{enumerate}
\item Every $S$-flat $R$-module is flat.
\item $S\subseteq U(R).$
\end{enumerate} 
\end{cor}
\proof $2.\Rightarrow 1.$ Obvious.

 $1.\Rightarrow 2.$ Let $s\in S$. By Theorem \ref{2thm-S-VNR}, $s$ is von Neumann regular, so there exists $r\in R$ such that $rs^2=s$. Since $s\in reg(R)$, $rs=1$. Hence, $s\in U(R)$. \cqfd 

 
 

\begin{exmp}\label{2exmp-S-VNR} Let $R_1$ be a von Neumann regular ring, $R_2$ be any commutative ring, $S_1$ be any multiplicative subset of $R_1$, and consider the ring $R=R_1\times R_2$ with the multiplicative subset $S = S_1\times U(R_2)$. Then every $S$-flat $R$-module is flat. Indeed, let $M=M_1\times M_2$ be an $S$-torsion $R$-module, $M_1$ (resp., $M_2$) is an $R_1$-module (resp., $R_2$-module). Then $M_2$ is $U(R_2)$-torsion, so $M_2=0$. Since $R_1$ is von Neumann regular, $M_1$ is a flat $R_1$-module. Hence, $M=M_1\times 0$ is a flat $R$-module (see \cite[Exercise 2.15, page 142.]{Wan1}). 
\end{exmp}

\section{S-perfect rings} 
It is  well-known  that a ring $R$ is perfect if and only if every flat $R$-module is projective \cite[Theorem 3.10.22]{Wan1}. So, it is natural to ask the following question: When every $S$-flat $R$-module is projective? This question leads us to introduce the following  notion.

\begin{defn} \label{3def-S-perfect} We say that $R$ is $S$-perfect if each $S$-flat $R$-module is a projective $R$-module.
\end{defn}

 Thus, since every flat $R$-module is an $S$-flat $R$-module, every $S$-perfect ring is perfect. Next, in Example \ref{3exmp-perfect-not-s-perfect}, we give an example of a perfect ring which is not $S$-perfect. It is  well-known  that if $R$ is a perfect ring, then $R_S$ is a perfect ring. Thus, if $R$ is an $S$-perfect ring, $R_S$ is a perfect ring. However, the converse is not necessarily the case as shown in Example \ref{3exmp-Rs-perfect-not-S-perfect}.  Recall that a ring $R$ is semisimple (resp., von Neumann regular) if every $R$-module is projective (resp., flat) \cite[page 212]{Rot1}. We know that  $R$ is semisimple if and only if $R$ is perfect and von Neumann regular. It is easy to verify that $R$ is semisimple if and only if $R$ is $S$-perfect and $R_S$ is von Neumann regular (Recall that $R_S$ is a von Neumann regular ring if and only if every $R$-module is $S$-flat).
 

\begin{exmp}\label{3exmp-Rs-perfect-not-S-perfect}
Let $R$ be a commutative von Neumann regular ring which is not semisimple. If $S=R \setminus p$ where $p$ is a prime ideal of $R$, then $R_S$ is a field (then a perfect ring) but $R$ not $S$-perfect because it is not perfect.
\end{exmp}


Recall that a ring $R$ is $S$-almost perfect if and only if every flat $R$-module is $S$-strongly flat \cite[Theorem 7.9]{Baz1}. It is clear that any $S$-perfect ring is $S$-weakly von Neumann regular and $S$-almost perfect. Moreover, if $R_S$ is projective, the converse is true.

\begin{thm}\label{3thm-S-perf-almostperf-vnr}    The following assertions are equivalent:
\begin{enumerate}
\item $R$ is $S$-perfect.
\item $R$ is $S$-weakly von Neumann regular and $S$-almost perfect, and $R_S$ is a projective $R$-module.
\item $R$ is perfect and  $S$-weakly von Neumann regular. 

 \end{enumerate}
\end{thm}

\proof 
$1.\Rightarrow 2.$ Follows immediately by the definitions ($R_S$ is projective because it is flat \cite[Theorem 4.80]{Rot1} and $R$ is perfect). 

$2.\Rightarrow 1.$ Follows from the definitions and Proposition \ref{2prop-proj-are-S-proj}.
 
$2.\Leftrightarrow 3.$  Follows from Corollary \ref{2cor-perfect-S-almost-perfect}.

\begin{cor}\label{3cor-s-perfect-and-vnr-elements} A ring $R$ is $S$-perfect if and only if $R$ is perfect and $R/Rs$ is a semisimple ring for any $s\in S$.
\end{cor}

\proof $(\Rightarrow )$ Let $s\in S$. $R/Rs$ is von Neumann regular by Theorem \ref{2thm-S-VNR} and perfect by \cite[Corollary 3.10.23]{Wan1}. Hence, it is semisimple. 


 $(\Leftarrow)$ This follows form Theorem \ref{3thm-S-perf-almostperf-vnr} (3) and Theorem \ref{2thm-S-VNR} (3). \cqfd

 It follows form Corollary \ref{3cor-s-perfect-and-vnr-elements} and Lemma \ref{2-recall-vnr} that a necessary condition for $R$ to be $S$-perfect is that $S\subseteq vnr^*(R)$. This fact allows us to construct a first example of a perfect ring which is not $S$-perfect.

\begin{exmp}\label{3exmp-perfect-not-s-perfect}Let $R = R_1\times R_2$ be the direct product of commutative perfect rings. Let $J(R_1)$ be the Jacobson radical  of $R_1$, $a_1 \in J(R_1)$, $a_2\in U(R_2)$ the set of units of $R_2$, and $a=(a_1,a_2)$. Let $S=\{a^n \mid n \geq 0\}$.  Then, $R$ is perfect, but it is not $S$-perfect.
\end{exmp}

\proof This follows from the following observations: 

1. It is clear that $S$ is a multiplicatively closed set, and $0\notin S$ because $a_2$ is assumed to be invertible, so $S$ is not trivial.

2. $R$ is  perfect by \cite[Theorem 3.10.22]{Wan1}.

3. It is easy to see that $a\in vnr(R)$ if and only if $a_1\in vnr(R_1)$ and $a_2\in vnr(R_2)$. Since $R_1$ is perfect, $J(R_1)$, the  Jacobson radical of $R_1$,  is $T$-nilpotent. Then, $a_1$ is a non-zero nilpotent, so $a_1 \notin vnr(R_1)$ \cite[Theorem 2.1, (3)]{And1}. Hence, $a\notin vnr(R)$. Thus, $S\nsubseteq vnr(R)$. Hence $R$ is not $S$-perfect. \cqfd

 Also, it is evident that every semisimple ring is $S$-perfect for any multiplicative set $S$. We use the following result to give an example of an $S$-perfect ring which is not semisimple,

\begin{prop}\label{3prop-s-perfectinproduct} Let $R = R_1 \times R_2\cdots\times R_n$ be the direct product of commutative rings $R_1,...,R_n$, and $S$ be a multiplicative subset of $R$. The following assertions are equivalent:
\begin{enumerate}
\item $R$ is $S$-perfect.
\item For each $i=1,...,n$, either $0\in e_iS$ and $R_i$ is a semisimple ring or $0\notin e_iS$ and $R_i$ is an $e_iS$-perfect ring. 
\end{enumerate}
\end{prop}

\proof This follows directly form Proposition \ref{2prop-locali-of-product} and Corollary \ref{2cor-S-projandS-flatinproduct}. We only sketch the proof. Assume that $R$ is $S$-perfect. Let $i\in \{1,... ,n\}$. Suppose $0\in e_i S$. Let $M_i$ be an $R_i$-module. By Proposition \ref{2prop-locali-of-product}, $(0\times M_i\times 0)_S=0$, where the $i$-th component of $0\times M_i\times 0$ is $M_i$ and the others are $0$. Then, $0\times M_i\times 0$ is an $S$-flat module. Since $R$ is $S$-prefect, $0\times M_i\times 0$ is a projective $R$-module, so $M_i$ is a projective $R_i$-module (see \cite[Exercise 2.15, page 142.]{Wan1}). Hence, $R_i$ is a semisimple ring. Now if $0\notin e_iS$, then $e_iS$ is multiplicative subset of $R_i$. Let $M_i$ be an $e_iS$-flat $R_i$-module, then, by corollary \ref{2cor-S-projandS-flatinproduct}, $0\times M_i\times 0$ is an $S$-flat $R$-module. Since $R$ is $S$-perfect, $0\times M_i\times 0$ is projective. Then, $M_i$ is a projective $R_i$-module. The converse holds by using the same arguments.    \cqfd

\begin{cor}\label{3cor-s-perfectinproductpratique} Let $ R = R_1 \times R_2$ be direct product of rings $R_1$ and $R_2$ and $S=S_1\times S_2$ be a cartesian product of multiplicative sets $S_i$ of $R_i$, $i=1,2$. Then, $R$ is $S$-perfect if and only if $R_i$ is $S_i$-perfect for every $i = 1, 2$.
\end{cor}


The promised $S$-perfect non-semisimple ring is given as follows.

\begin{exmp}\label{3exmp-s-peferfectnonsemisimple}
 Let $R_1$ be a semisimple ring and $R_2$ be a non-semisimple perfect  ring. Then $R = R_1 \times R_2$ is not a semisimple ring. Let $S$ be any multiplicative subset of $R_1$. Set $S = S\times \{1\}$
which is a multiplicative subset of $R$. Then $R$ is an $S$-perfect ring by Corollary \ref{3cor-s-perfectinproductpratique}, but it is not semisimple.
\end{exmp}

 We summarize some results given above:

\begin{center}
$R$ is semisimple $\Rightarrow$  $R$ is $S$-perfect  $\Rightarrow$ $R$ is perfect rings $\Rightarrow$ $R_S$ is perfect rings 
 \end{center}
 
Now, it is not difficult to see that all the above implications are irreversible. We then return to complete the study of $S$-perfect rings. We have the following homological characterization,

\begin{thm}\label{3thm-homologicalchara-sperfect}    The following assertions are equivalent:
\begin{enumerate}
\item $R$ is $S$-perfect.
\item Every $S$-flat module is $S$-projective and projective.
\item $R_S$ is a perfect ring, projective (cyclic) as $R$-module and every $S$-projective module is projective.
\item $R_S$ is a perfect ring and, for every $R$-module $M,$ if $M_S$ is $R_S$-projective, then $M$ is a projective $ R$-module.
\item $\pd_R(M)=\fd_{R}(M_S)$ for any $R$-module $M.$
 \end{enumerate}
\end{thm}

\proof $1.\Rightarrow 2.$ Let $M$ be an $S$-flat module, then $M$ is projective by Definition \ref{3def-S-perfect}. By Proposition \ref{2-prop-S-flat-charac} $M_S$ is a flat $R$-module, then $M_S$ is a projective $R$-module because every $S$-prefect ring is perfect.
 
 $2.\Rightarrow 1.$ This is obvious. 
 
$1.\Rightarrow 3.$ It is clear that $R_S$ is a perfect ring (see the discussion after Definition \ref{3def-S-perfect}). $R_S$ is  projective because $R_S$ is a flat $R$-module and every flat $R$-module is $S$-flat. Since $R$ is $S$-perfect, it is $S$-weakly von Neumann regular. Then, by Theorem \ref{2thm-S-VNR}, the canonical homomorphism $\pi : R \to R_S$ is an epimorphism. So $R_S$ is a cyclic $R$-module. The last assertion follows from the fact that every $S$-projective $R$-module is an $S$-flat $R$-module.

 $3.\Rightarrow 4.$ Let $M$ be an $R$-module. Assume that $M_S$ is a projective $R_S$-module. Since $R_S$ is a projective $R$-module, $M_S$ is a projective $R$-module, i.e., $M$ is $S$-projective. Hence $M$ is projective by hypothesis. 

$4.\Rightarrow 1.$ Let $M$ be an $S$-flat $R$-module, then $M_S$ is a flat $R_S$-module. Since $R_S$ is perfect, $M_S$ is a projective $R_S$-module, so $M$ is a projective $R$-module, as desired. 

$5.\Rightarrow 1.$  Follows form Proposition \ref{2-prop-S-flat-charac}.

$1.\Rightarrow 5.$ If $\pd_R(M)=\infty$, so is $\fd_R(M_S)$. Indeed, assume $\fd_R(M_S)=n <\infty$. If $\cdots\to P^1 \to P^0 \to M\to 0$ is a projective resolution of $M$, then, since $R_S$ is a flat $R$-module,
$$0 \to K^{n-1}_S\to P^{n-1}_S \to\cdots\to  P^1_S \to P^0_S\to M_S\to 0$$

is a  resolution of $M_S$ by $R$-modules, where $K^{n-1}$ is the ($n - 1)$st syzygy. Since $P^i_S$ are flat $R$-modules and $\fd_R(M_S)=n$, $K^{n-1}_S$ is a flat $R$-module. Then, by proposition \ref{2-prop-S-flat-charac}, $K^{n-1}$ is an $S$-flat $R$-module, so it is projective because $R$ is $S$-perfect. Then 

$$0 \to K^{n-1}\to P^{n-1} \to\cdots\to  P^1 \to P^0\to M\to 0,$$
  
  is a finite projective resolution of $M$, it is a contradiction.
  
 If $\pd_R(M)=n <\infty$. Consider a projective resolution of $M$ of length $n$ :
  $$0 \to K^{n-1}\to P^{n-1} \to\cdots\to  P^1 \to P^0\to M\to 0,$$
then, 
$$0 \to K^{n-1}_S\to P^{n-1}_S \to\cdots\to  P^1_S \to P^0_S\to M_S\to 0,  \,\,\,\, (1)$$
 is a flat resolution of $M_S$ of length $n$. So, $\fd_R(M_S)\leq n$. If  $\fd_R(M_S)=m< n$, then $K^{m-1}_S$, the ($m - 1)$st   syzygy of $(1)$, is a flat $R$-module, i.e., $K^{m-1}$ is $S$-flat. Since $R$ is $S$-perfect, $K^{m-1}$ is a projective $R$-module. Then  $\pd_R(M)\leq m< n$, a contradiction. Thus $\fd_R(M_S)= n=\pd_R(M)$, as desired.
\cqfd   

It is well-known  that if $R$ is a perfect ring, then $R_S$ is a perfect ring and $R_S$ is a projective $R$-module (recall that  $R_S$ is a flat $R$-module). Thus, Theorem \ref{3thm-homologicalchara-sperfect} implies the following characterization of $S$-perfectness over perfect rings.

\begin{cor}\label{3cor-homologicalchara-sperfect}  If $R$ is a perfect ring, then the following assertions are equivalent:
\begin{enumerate}
\item $R$ is $S$-perfect, 
\item Every $S$-projective module is projective,
\item $\pd_R(M)=\pd_{R}(M_S)$ for any $R$-module $M.$
 \end{enumerate}
\end{cor}

\begin{cor}\label{cor-316} If $R$ is an $S$-perfect ring, then $$\gldim(R)=\wdim(R)=\gldim(R_S)=\wdim(R_S)$$
\end{cor}

\proof Since $R$ is $S$-perfect, $R$ and $R_S$ are perfect rings. Then, $\wdim(R)=\gldim(R)$ and  $\wdim(R_S)=\gldim(R_S)$. Also, the inequalities $\wdim(R_S)\leq \wdim(R)\leq \gldim(R)$ are always verified \cite[Proposition 8.21 and Proposition 8.23.]{Rot1}. We just need to prove the reverse inequality of $\wdim(R_S)\leq \gldim(R)$, we may assume that $\wdim(R_S) <\infty$. Then $\wdim(R_S)=\gldim(R_S)=0$ \cite[Theorem 3.10.25]{Wan1}. Let $M$ be an $R$-module. Since $R$ is $S$-perfect, by Theorem \ref{3thm-homologicalchara-sperfect}, $\pd_ R(M)=\fd_R(M_S)= 0$. Hence, $\gldim(R)= \wdim(R_S)=0$.\cqfd

  We say that an $R$-module $M$ is $S$-divisible if $M=sM$ for every $s\in S$ \cite[Definition 1.6.10]{Wan1}.

\begin{lem}\label{3lem-s-divisivle}
\begin{enumerate}
\item Let $R = R_1 \times R_2 \times \cdots \times R_n$ be a direct product of rings and let $M=M_1\times \cdots\times M_n$ be an $R$-module. Then $M$ is $S$-divisible if and only if each $M_i$ is $e_iS$-divisible.
\item If $0\in S$, then an $R$-module $M$ is $S$-divisible if and only if $M=0$. 
\item Let $S\subseteq S'$ be multiplicative subsets of $R$ and $M$ an $R$-module. If $M$ is $S'$-divisible, then it is $S$-divisible.
\end{enumerate}
\end{lem}

\proof Straightforward.

\begin{lem}\label{3lem-perfect-vonelements}
Let $J$ be the Jacobson radical of $R$. Then the following statements holds:
 \begin{enumerate}
 \item If $R$ is perfect, then $0\in S$ if and only if $J\cap S \neq  \emptyset.$
 \item If $J$ is prime, then $vnr^*(R)=U(R)$.  
 \end{enumerate}
\end{lem}
\proof 1. $(\Rightarrow)$ Clear. $(\Leftarrow)$ Let $x\in J\cap S$, since $R$ is perfect, $J$ is $T$-nilpotent. Then there is $n\in \N$ such that $x^n=0$, hence $0\in S$.

2. Let $x\in vnr(R)$, then $x(1-ux)=0$ for some $u\in U(R)$ \cite[Theorem 2.2.]{And1}. If $J$ is prime then $x\in J$ or $1-ux\in J$. If $x\in J$, so is $ux$. Then  $1-ux \in U(R)$, so $x=0$. If $1-ux\in J$, then $ux\in U(R)$. Since $u\in U(R)$, $x\in U(R)$.  \cqfd

Recall that a ring $R$ is perfect if and only if $R = R_1 \times R_2 \times\cdots \times R_p$, with each $R_i$ is local perfect ring. We have the following characterization of $S$-perfect rings.

\begin{thm}\label{3thm-S-perfect-jacobson} The following assertions are equivalent:
\begin{enumerate}
\item $R$ is $S$-perfect.
\item $R$ is perfect and the Jacobson radical $J(R)$ is $S$-divisible.
\item $R = K_1\times \cdots\times  K_n\times R_1\times\cdots\times R_m$ such that $K_i$ are fields, $R_i$ are  perfect local rings, $0\in e_iS $ for each $i=1,..., n$ and $e_iS \subseteq U(R_i)$ for each $i=1,..., m$.  
\end{enumerate} 
 \end{thm} 

\proof $1.\Rightarrow 3.$ Assume that $R$ is an $S$-perfect ring. Then it is prefect, so, by \cite[Theorem 3.10.22]{Wan1}, $R = R_1 \times R_2 \times\cdots \times R_p$, where each $R_i$ is a local perfect ring. We may assume that there is $n< p$ such that $0\in e_i S \Leftrightarrow i\leq n $. Since $R = R_1 \times R_2 \times\cdots \times R_p$ is an $S$-perfect ring, by Proposition \ref{3prop-s-perfectinproduct}, $R_i$ is a semisimple ring for each $i\leq n$ and $R_i$ is an $e_iS$-perfect ring for each $i>n$. Since $R_i$ is local, $R_i$ is a field for each $i=1,...,n$. The last assertion follows from Corollary \ref{3cor-s-perfect-and-vnr-elements} and Lemma \ref{3lem-perfect-vonelements} (2). 

$3.\Rightarrow 1.$ This follows form Proposition \ref{3prop-s-perfectinproduct}.

 $2.\Rightarrow 3.$ Since $R$ is perfect, by \cite[Theorem 3.10.22]{Wan1}, $R = R_1 \times R_2 \times\cdots \times R_p$, where each $R_i$ is a local perfect ring. Also, we may assume that there is $n< p$ such that $0\in e_i S \Leftrightarrow i\leq n $. By \cite[Exercise 1.12., page 64]{Wan1}, $J(R)= J(R_1)\times \cdots\times J(R_p)$. Since $J(R)$ is $S$-divisible, by Lemma \ref{3lem-s-divisivle} (1), $J(R_i)$ is $e_iS$-divisible for each $i=1,..., p$, but, since $0\in e_iS$, $ J(R_i)=0$ for each $i\leq n$ by Lemma \ref{3lem-s-divisivle} (2). Thus, $R_i$ is a filed for each $i\leq n$. Let $i>n$, then $0\notin e_iS$. Since $R_i$ is perfect, by Lemma \ref{3lem-perfect-vonelements} (1), $J(R_i)\cap e_iS= \emptyset$. So, $e_iS\subseteq U(R_i)$ because $R_i$ is local.
 
 $3.\Rightarrow 2.$ $R$ is perfect because it is a direct product of perfect rings  \cite[Theorem 3.10.22]{Wan1}.  To see that $J(R)$ is $S$-divisible, note that $S\subseteq S'=K_1\times \cdots\times K_n\times U(R_1)\times\cdots\times U(R_m)$, $J(R)= 0\times J(R_1)\times\cdots\times J(R_m)$ and  each $J(R_i)$ is $e_iS'$-divisible. Therefore Lemma \ref{3lem-s-divisivle} $(3)$ gives the result. \cqfd

 \begin{cor}
 Let $ S' \subseteq S$ be multiplicative subsets of a commutative ring $R$. If $R$ is $S$-perfect, then it is $S'$-perfect.
 \end{cor}

\proof  Follows by Theorem \ref{3thm-S-perfect-jacobson} (2) and Lemma \ref{3lem-s-divisivle} (3). \cqfd
 

\begin{rem}\label{cor-322}
Let $R$ be a commutative non-semisimple perfect ring. Then, by \cite[Theorem 3.10.22]{Wan1}, $R = R_1 \times R_2 \times\cdots \times R_p$, where each $R_i$ is a local perfect ring. Then, either $R_i$ is not a field for each $i$, or there exists a non-trivial multiplicative subset $S$ of $R$ such that $R$ is an $S$-perfect ring (by "non-trivial" multiplicative subset $S$ we mean $0\notin S$ and $S\nsubseteq U(R)).$ Indeed, suppose that there exists $i$ such that $R_i$ is a field. We may assume that there exists $1 \leq m<n$ such that $R_i$ is field if and only if $i\leq m$. Let $S=U(R_1)\times \cdots \times U(R_m)\times 0$.  Then $R$ is $S$-perfect by Theorem \ref{3thm-S-perfect-jacobson}.
\end{rem}

Here, an other non-trivial example of a perfect (resp., an $S$-perfect) ring which is not $S$-perfect (resp., semisimple) with respect to a non-trivial multiplicative subset $S$ of $R$: 

\begin{exmp}\label{exmp-323} Let $K$ be a field and $R$ a non-semisimple local perfect ring, then $K\times R$ is a non-semisimple  perfect ring and we have,
\begin{enumerate}
\item $K\times R$ is not $(K^*\times 0)$-perfect
\item $K\times R$ is $(0\times U(R))$-perfect
\end{enumerate}
\end{exmp}

Next, we examine the transfer of the $S$-perfect property along a ring homomorphism.  Let $f : R \to R'$ be a ring homomorphism and $S$ a multiplicative subset of $R$. Then it is easy to see that $f(S)$ is a multiplicative subset of $R'$ when $ 0 \notin f(S)$. Recall that every $R'$-module $M$ acquires an $R$-module structure via the module action $rm= f(r)m$ \cite[Proposition 8.33.]{Rot1}. In particular, if $ R \hookrightarrow R'$ is a ring extension (that is, $R$ is a subring of $R'$) and $S$ is a multiplicative subset of $R$, then $S$ is also a multiplicative subset of $R'$.

\begin{prop}\label{3prop-S-pefect-in-homomorphism}
Let $ R  \hookrightarrow R'$ be a ring  extension such that $R'$ is finitely generated flat $R$-module and let $S$ be a multiplicative subset of $R$ (so of $R'$). If $R'$ is $S$-perfect, then so is $R$.
\end{prop}

Before proving Proposition \ref{3prop-S-pefect-in-homomorphism}, we establish the following lemma:

\begin{lem}\label{lem-325}
Let $f : R \to  R'$ be a ring homomorphism 	and $S$ a multiplicative subset of $R$  with $0 \notin  f(S)$. Let $M$ be an $R'$-module. Then there is an isomorphism of $R$-module $M_S\cong M_{f(S)}$.  
\end{lem}
\proof Recall that $M$ acquires an $R$-module structure via the module action $rm= f(r)m$. Define $\phi: M_S \to M_{f(S)}$ by $\phi(\frac{m}{s}) = \frac{m}{f(s)}$. Then $\phi$ is well-defined and injective because of the $R$-module structure on $M$. Also, the verification of the surjectivity is easy.  \cqfd

Now we are able to prove Proposition \ref{3prop-S-pefect-in-homomorphism}.\\
\textbf{Proof of Proposition \ref{3prop-S-pefect-in-homomorphism}}  Let $M$ be an $S$-flat $R$-module; that is, $M_{S}$ is a flat $R$-module. We show that $M$ is a projective $R$-module and then the result follows from Definition \ref{3def-S-perfect}. Since $R'$ is a flat $R$-module, $M_S\otimes_R R'$ is a flat $R'$-module. But, $$M_S \otimes_R R'\cong (R_S\otimes_R M)\otimes_R R' \cong R_S\otimes_R (M\otimes_R R'),$$ and, by Lemma \ref{lem-325}, $R_S\otimes_R(M\otimes_R R') \cong R'_S\otimes_{R'}(M\otimes_R R')$. so $M\otimes_R R'$ is an $S$-flat $R'$-module by Proposition \ref{2-prop-S-flat-charac}.
 Since $R'$ is $S$-perfect, $M\otimes_R R'$ is a projective $R'$-module.  Then, by \cite[Corollary 3.3]{Sha1}, $M$ is a projective $R$-module. \cqfd

\begin{cor}\label{3cor-S-pefect-in-isomomorphism}
Let $f : R \to  R'$ be a ring isomorphism and $S$ a multiplicative subset of $R$. Then, $R$ is an $S$-perfect ring if and only if $R'$ is $f(S)$-perfect.
\end{cor}

 Let $R\ltimes M$ be the trivial ring extension of $R$ by an $R$-module $M$. If $S$ is a multiplicative subset of $R$ and $N$ is a submodule of $M$, then $S\ltimes N$ is a multiplicative subset of $R\ltimes M$. Conversely, if $S$ is a multiplicative subset of $R\ltimes M$, then $S_1=\{r\in R\mid (r,m)\in  S \text{ for some $m$ in $M$} \}$ is a multiplicative subset of $R$. Indeed, it is clear that $1\in S_1$ and $S_1$ is closed under multiplication, so we only need to check  that $0\notin S_1$, but, if $0\in S_1$, then $(0,m)\in S$ for some $m\in M$ and so $(0,0)=(0,m)^2\in S$, a contradiction.

\begin{lem}\label{3lemoftrivial-S-pefect-in-homomorphism}
 Let $S$ be a multiplicative subset of $R\ltimes M$. If $R\ltimes M$ is an $S$-perfect ring, then $R$ is an $S_1$-perfect, where $S_1=\{r\in R\mid (r,m)\in S \text{ for some m in M}\}$. In particular, if $S$ is a multiplicative subset of $R$ and $N$ is a submodule  of $M$ such that $R\ltimes M$ is an $S\ltimes N$-perfect ring, then $R$ is an $S$-perfect ring.
\end{lem}

 \proof Since every $S$-perfect ring is perfect, $R\ltimes M$ is perfect. Then, by \cite[Proposition 1.15., page 23]{Fos1}, $R$ is perfect. So, it suffices by Theorem \ref{3thm-S-perfect-jacobson} to show that $J(R)$ is $S_1$-divisible. Let $s\in S_1$ and $r\in J(R)$. We have $(s,m)\in S$ for some $m\in M$ and $(r,0)\in J(R\ltimes M)= J(R)\ltimes M$. Since $R\ltimes M$ is an $S$-perfect, by Theorem \ref{3thm-S-perfect-jacobson}, $J(R\ltimes M)$ is $S$-divisible. Then, there exists $(r',m') \in J(R\ltimes M)$  such that $(r,0)=(s,m)(r',m')$. Hence, $r=sr'$ with $r'\in J(R)$. Since $s$ and $r$ are arbitrary in $S_1$ and $J(R)$, respectively, $J(R)$ is $S_1$-divisible.  \cqfd

  Next result shows that, unlike the "perfect" property \cite[Proposition 1.15., page 23]{Fos1}, the "$S$-perfect" property depends on the structure of the $R$-module $M$ considered in the trivial extension $R\ltimes M$. 

\begin{thm}\label{3thm-s-perfect-in-trival-extention}
 Let $S$ be a multiplicative subset of $R\ltimes M$.  Then, the following statements are equivalent: 
 \begin{enumerate}
 \item $R\ltimes M$ is an $S$-perfect ring.
 \item $R$ is an $S_1$-perfect ring and $M$ is an $S_1$-divisible $R$-module, where $S_1=\{r\in R\mid (r,m)\in S \text{ for some m in M}\}$. 
 \end{enumerate}  
\end{thm}
  
\proof $1.\Rightarrow 2.$ Assume that $R\ltimes M$ is $S$-perfect, then, by Lemma \ref{3lemoftrivial-S-pefect-in-homomorphism}, $R$ is $S_1$-perfect. Then, by Theorem \ref{3thm-S-perfect-jacobson}, $R = K_1\times \cdots\times  K_n\times R_{n+1} \times\cdots\times R_{m}$ such that $K_i$ are fields, $R_i$ are  perfect local rings, $0\in e_iS_1$ for each $i=1,..., n$ and $e_iS_1 \subseteq U(R_i)$ for each $i=n+1,..., m$, where $e_i = (0, ..., 1, ..., 0)$, the $i$-$th$ component of $e_i$ is $1$ and the others are $0$. Write  $M= E_1\times \cdots\times  E_n\times M_{n+1}\times\cdots\times M_{m}$. Clearly, $R\ltimes M=(K_1\ltimes E_1)\times \cdots\times  (K_n\ltimes E_n)\times (R_{n+1}\ltimes M_{n+1})\times\cdots\times (R_{m}\ltimes M_{m})$. Let $e'_i=((0,0),...,(1,0),...,(0,0))$; that is, the $i$-$th$ component of $e'_i$ is $(1,0)$ and the others are $0$. Now, $(0,0)\in e'_i S$ for any $i=1,...,n$. Indeed, for $i\in \{1,...,n \}$, we have $0\in e_iS_1$, then $(0,n)\in e'_iS$ for some $n\in e_i M$. Hence, $(0,0)=(0,n)^2\in e'_iS$. Since $R\ltimes M$ is $S$-perfect, by Proposition \ref{3prop-s-perfectinproduct}, $K_i\ltimes E_i$ are semisimple rings. Hence, $E_i=0$ for any $i=1,...,n$. Thus, in light of lemma \ref{3lem-s-divisivle}, $M$ is $S_1$-divisible if and only if $M_i=e_i M$ is $e_iS_1$-divisible for any $i=n+1,...,m$.  But, since $e_iS_1 \subseteq U(R_i)$ for each $i=n+1,..., m$, $M_i=e_i M$ is $e_iS_1$-divisible for any $i=n+1,...,m$. Hence, $M$ is $S_1$-divisible.

$2.\Rightarrow 1.$ Since $R$ is $S_1$-perfect, $R$ is perfect and then $R\ltimes M$ is a perfect ring by \cite[Proposition 1.15., page 23]{Fos1}. Then, in light of Theorem \ref{3thm-S-perfect-jacobson}, we need just to prove that $J(R\ltimes M)= J(R)\ltimes M$ is $S$-divisible. Let $(s,n)\in S $ and $(r,m)\in J(R)\ltimes M$. We have to find $(r',m')\in J(R)\ltimes M$ satisfied $(r,m)=(s,n)(r',m')$. Since $J(R)$ and $M$ are  $S_1$-divisible, there exist $r'\in J(R)$ and $m'\in M$ such that $r=sr'$ and $m=sm'$. By Corollary \ref{3cor-s-perfect-and-vnr-elements}, $s$ is a von Neumann regular, so, $s=us^2$ for some $u\in R$. Finally, $(r,m)=(s,n)(usr',-ur'n+usm')$, as desired. \cqfd

\begin{cor}\label{3cor-s-perfect-in-trival-extention}
 Let $S$ be a multiplicative subset of $R$ and $N$ is a submodule  of $M$.  Then, the following statements are equivalent: 
 \begin{enumerate}
 \item $R\ltimes M$ is an $S\ltimes N$-perfect ring.
 \item $R\ltimes M$ is an $S\ltimes 0$-perfect ring.
 \item $R$ is an $S$-perfect ring and $M$ is an $S$-divisible $R$-module.
\end{enumerate}   
\end{cor}

\begin{exmp}\label{2exmp-S-VNR} Let $R$ be a commutative semisimple ring and $M$ be a non-zero $R$-module. Consider the ring $R'=R\times R $ with the multiplicative subset $S' = 0\times U(R)$. Then $R'$ is an $S'$-perfect ring, but, $R'\ltimes M'$ is not $S'\ltimes 0$-perfect, where $M'=M\times M$. Indeed, $R'$ is $S'$-perfect because it is semisimple. Since, by Lemma \ref{3lem-s-divisivle} (1) and (2), $M'$ is not $S'$-divisible, $R'\ltimes M'$ is not $S'\ltimes 0$-perfect by Corollary \ref{3cor-s-perfect-in-trival-extention}. 
\end{exmp}

Driss Bennis:  Faculty of Sciences, Mohammed V University in Rabat, Rabat, Morocco.

\noindent e-mail address: driss.bennis@um5.ac.ma; driss$\_$bennis@hotmail.com

Ayoub Bouziri: Faculty of Sciences, Mohammed V University in Rabat, Rabat, Morocco.

\noindent e-mail address: ayoub$\_$bouziri@um5.ac.ma


\begin{thebibliography}{999}

\bibitem{Abu1} E. Abu Osba, M. Henrikson and O. Alkam,\textit{ Combining local and von Neumann regular rings}, Commun. Algebra \textbf{32} (2004), 2639-2653.

\bibitem{And1} D. F. Anderson and A. Badawi, \textit{Von Neumann regular and related elements in commutative rings}, Algebra Colloq. \textbf{19} (2012), 1017-1040.

\bibitem{And2} D. D. Anderson and T. Dumitrescu,\textit{ S-Noetherian rings}, Commun. Algebra \textbf{30} (2002), 4407-4416.

\bibitem{Ati1} M. F. Atiyah and I.G. Macdonald, \textit{Introduction to Commutative Algebra}, Addison-Wesley Pub-
lishing Company, London, (1969).

\bibitem{Baz1}  S. Bazzoni and L. Positselski, \textit{S-almost perfect commutative rings}, J. Algebra \textbf{532} (2019),
323-356.
\bibitem{Benn1} D. Bennis and M. El Hajoui, \textit{On S-coherence}, J. Korean Math. Soc. \textbf{55} (2018), 1499-1512.
 
\bibitem{Cha1} S. U. Chase, \textit{Direct products of modules}, Trans. Amer. Math. Soc. \textbf{97} (1960), 457-473. 

\bibitem{Fos1} R. M. Fossum, P. Griffith and I. Reiten, \textit{ Trivial Extensions of Abelian Categories}, Lecture Notes \textbf{456}, Springer-Verlag, New York (1975).

\bibitem{Gla1} S. Glaz, \textit{Commutative Coherent Rings}, Lecture Notes in Mathematics, \textbf{1371}, Spring-Verlag, Berlin (1989).
   



\bibitem{Pal1} I. Palmer and J. E. Roos, \textit{Explicit formulae for the global homological dimensions of trivial extensions of rings}, J. Algebra \textbf{27} (1974), 380-413.

\bibitem{Pos1}  L. Positselski and A. Slávik, \textit{On strongly flat and weakly cotorsion modules}, Math. Z. \textbf{291} (3-4) (2019), 831-875.

\bibitem{Qi1} W. Qi, X. Zhang and W. Zhao, \textit{New Characterizations of S-coherent rings}, J. Algebra Appl. \textbf{22} (2023), 2350078 (14 pages).



\bibitem{Rot1} J. Rotman, \textit{An Introduction to Homological Algebra}, Academic Press, New York, (2009).

\bibitem{Sha1} A. Shamsuddin, \textit{Finite normalizing extensions}. J. Algebra, \textbf{151}, (1992), 218-220.
\bibitem{Wan1} F. G. Wang and H. Kim, \textit{Foundations of Commutative Rings and Their Modules}, Singapore, Springer, (2016). 

\end{thebibliography}
\end{document}